\documentclass[reqno]{amsart}

\newcommand{\field}[1]{\ensuremath{{\mathbb #1}}}
\newcommand{\Quat}{\field{H}}
\newcommand{\R}{\field{R}}
\newcommand{\C}{\field{C}}
\newcommand{\N}{\field{N}}
\newcommand{\F}{\field{F}}
\newcommand{\Q}{\field{Q}}
\newcommand{\Proj}{\ensuremath{{\mathbb P}}}
\newcommand{\fpn}{\ensuremath{\F\Proj^n}}
\newcommand{\fp}{\ensuremath{\F\Proj}}

\newcommand{\eqskip}{\quad}

\newcommand{\cD}{\ensuremath{{\mathcal D}}}
\newcommand{\cM}{\ensuremath{{\mathcal M}}}
\newcommand{\cL}{\ensuremath{{\mathcal L}}}

\newcommand{\Pol}{\ensuremath{\mathrm{Pol}}}
\newcommand{\polfd}{\ensuremath{\Pol_\F(d)}}

\newcommand{\dm}{\,d}

\newcommand{\abs}[1]{\ensuremath{\lvert#1\rvert}}

\DeclareMathOperator{\id}{id}

\DeclareMathOperator{\rank}{rk}
\DeclareMathOperator{\tr}{tr}
\DeclareMathOperator{\res}{res}

\newtheorem{theorem}{Theorem}[section] 
\newtheorem{thm}{Theorem} 
\newtheorem{lemma}[theorem]{Lemma}
\newtheorem{prop}[theorem]{Proposition}
\newtheorem{cor}[theorem]{Corollary}

\theoremstyle{remark}
\newtheorem{remark}[theorem]{Remark}
\newtheorem{example}[theorem]{Example}

\numberwithin{equation}{section}

\begin{document}


\author{Yu. I. Lyubich}
\title{On tight projective designs}
 
\begin{abstract}
  It is shown that among all tight designs in $\fpn\neq\R\Proj^1$,
  where $\F$ is $\R$ or $\C$, or $\Quat$ (quaternions), only 5-designs
  in $\C\Proj^1$ \cite{lyushatdedicata} have irrational angle
  set. This is the only case of equal ranks of the first and the
  last irreducible idempotent in the corresponding Bose-Mesner algebra.


\end{abstract}

\maketitle

\section{Introduction}
\label{sec:intro}

A well known theorem of Bannai and Hoggar \cite{bannaihoggarsqfree}
states that there are no tight $t$-designs in $\fpn\neq\R\Proj^1$ if
$t\geq 6$. Moreover, a theorem of Hoggar \cite{hoggar45designs} states
the same for $t\geq 4$ if $\F\neq\R$. Surprisingly, a tight 5-design
in $\C\Proj^1$ has been constructed in \cite{lyushatdedicata}, so
Hoggar's theorem has to be corrected. The results of
\cite{bannaihoggarsqfree} and \cite{hoggar45designs} are essentially
based on Theorem 2.6(c) \cite{hoggaroctonions} that states that the
angle set of every tight $t$-design in $\fpn\neq\R\Proj^1$ is
rational. But it is not rational for the $5$-design constructed in
\cite{lyushatdedicata}.

In the present paper we investigate this contradiction and prove that
{\em the only cases where the angle set is not rational are} 
\begin{enumerate}
  \item $\F=\C$, $n=1$, $t=5$ {\em and} 
  \item $\F=\R$, $n=1$, $t\neq   1,2,3,5$.
\end{enumerate}
A fortiori, there are no complications in \cite{bannaihoggarsqfree}
where $t\geq 6$ by assumption.

Our principal observation is that if $t=2s-1$, $s\geq 2$ then the last
irreducible idempotent $L_s$ in the corresponding Bose-Mesner algebra is not
$E_s$ from the proof of Theorem 2.6(c) \cite{hoggaroctonions}
(actually, from \cite{neumaier}). Nevertheless, $\rank L_s \neq \rank
E_1$, except for our case (1). This ``critical inequality'' implies the
rationality of the angle set, similarly to the argument in 
\cite{hoggaroctonions}. This material is concentrated in Section
\ref{sec:4critineq} of the present paper, while Sections
\ref{sec:2projtdesigns} and \ref{sec:3bosemesner} contain all the necessary
background and preliminary analysis.

\section{Projective $t$-designs}
\label{sec:2projtdesigns}

For the reader's convenience we basically use the same notation as in
\cite{hoggardesignsproj} and other related papers. Let us recall this notation.
In particular, let
\begin{equation*}
  \F\in\{\R,\C,\Quat\}; \eqskip m=\frac{1}{2}(\F:\R) = 
  \begin{cases}
    1/2 & \F=\R \\
    1   & \F=\C \\
    2   & \F=\Quat
  \end{cases}; \eqskip N=m(n+1).
\end{equation*}
The number $2N$ is nothing but the real (topological) dimension of the
$\F$-linear space $\F^{n+1}$. The latter consists of all $(n+1)\times
1$ matrices (columns) over $\F$ with the standard addition and
multiplication by scalars $\tau\in\F$ from the right (for
definiteness). As usual, the inner product of $a,b\in\F^{n+1}$ is
$a^*b$ where $a^*$ is the row conjugate transpose to $a$. Accordingly,
the set
\begin{equation*}
  S^{2N-1} = \left\{ a : a^*a = 1\right\}
\end{equation*}
is the unit sphere in $\F^{n+1}$. A quotient set of the sphere with
respect to the equivalence relation $a_1\sim a_2 \Longleftrightarrow
a_1=a_2\lambda$, $\lambda\in\F$, $\abs{\lambda}=1$, is the projective
space $\fpn$. The ``inner product'' $(\hat a, \hat b) = \abs{a^*b}^2$
in $\fpn$ is well-defined through the natural mapping $a\mapsto\hat
a$ from $S^{2N-1}$ onto $\fpn$. Obviously, $(\hat b, \hat a) = (\hat
a, \hat b)$ and $0 \leq (\hat a, \hat b) \leq 1$ with the equality
$(\hat a, \hat b)=1$ if and only if $\hat a = \hat b$. For every nonempty
$X\subset\fpn$ its angle set is
\begin{equation*}
  A(X) = \left\{(x,y) : x,y\in X,x\neq y\right\}
\end{equation*}
The related combinatorial parameters are
\begin{equation*}
  s = \abs{A(X)},\eqskip e=\abs{A(X)\setminus\{0\}},\eqskip \epsilon=s-e = \abs{A(X)\cap\{0\}}.
\end{equation*}

Let $P_i^{(\alpha,\beta)}(\tau)$ be the Jacobi polynomials
\cite{szegogabor} such that
\begin{equation}
  \label{eq:2.0jacpoly}
  \deg P_i^{(\alpha,\beta)} = i, \eqskip P_i^{(\alpha,\beta)}(1) = \frac{(\alpha+1)_i}{i!}
\end{equation}
where
\begin{equation*}
  (\alpha+1)_i = \prod_{l=1}^i(\alpha+l),\eqskip (\alpha+1)_0 = 1.
\end{equation*}
In particular, $P_0^{(\alpha,\beta)}(\tau) \equiv 1$. In what follows
we fix 
\begin{equation}
  \label{eq:2.0palphabeta}
  \alpha=N-m-1,\eqskip\beta=m-1
\end{equation}
and set
\begin{equation}
  \label{eq:2.0pp}
  P_i(\xi) = P_i^{(\alpha,\beta)}(2\xi-1),
\end{equation}
for short. A finite nonempty subset $X\subset\fpn$ is called a
{\em $t$-design} if
\begin{equation}
  \label{eq:2.1tdesign}
  \sum_{x\in X} P_i((x,y)) = 0, \eqskip y\in X, \eqskip 1\leq i \leq t.
\end{equation}

Let $X$ be a $t$-design and let
\begin{equation*}
  R_e^\epsilon(\xi) = \frac{(N)_s}{(m)_s}P_e^{(\alpha+1,\beta+\epsilon)}(2\xi-1).
\end{equation*}
In particular,
\begin{equation}
  \label{eq:2.1p}
  R_e^\epsilon(1) = \frac{(N)_s(N-m+1)_e}{(m)_s e!}
\end{equation}
The following theorems are fundamental, see \cite{bannaiextremal},
\cite{bannaihoggarsymmetric},
\cite{hoggardelsartespaces}. (Cf. \cite{delsartegoethalsseidel} for
the spherical designs.)

\begin{thm}
  The inequalities
  \[
    t\leq s+e,\eqskip \abs{X}\geq R_e^\epsilon(1)
  \]
  hold, and the equalities
  \[
    t=s+e,\eqskip \abs{X} = R_e^\epsilon(1)
  \]
  are equivalent.
\end{thm}

In the latter case the $t$-design $X$ is called $tight$. Note that
$t=s+e$ is equivalent to $e=[t/2]$, $\epsilon=\res_2(t)$.

\begin{thm}
  If $X$ is a tight $t$-design then $A(X)$ coincides with the set of
  roots of the polynomial $\xi^\epsilon R_e^\epsilon(\xi)$.
\end{thm}

Recall that these roots are simple and lie on $(0,1)$.


\begin{thm}
  Let $X$ be a subset of $\fpn$ such that $\abs{X}=R_e^\epsilon(1)$ and
  $A(X)$ coincides with the set of roots of $\xi^\epsilon
  R_e^\epsilon(\xi)$, then $X$ is a tight $(2e+\epsilon)$-design.
\end{thm}

The projective $t$-designs can be characterized as the averaging sets
in the sense of \cite{seymourzaslavsky} for suitable spaces of
functions on $\fpn$. Usually, these spaces are described in
terms of harmonic analysis but we prefer a more elementary approach
\cite{lyushatquaternions}, \cite{lyushatpolyfunc}. 

We say that a mapping $\phi:\fpn\rightarrow\C$ is a {\em
  polynomial function} if it is of the form
\begin{equation*}
  \phi(\hat a) = \psi(a),\eqskip a\in S^{2N-1},
\end{equation*}
where $\psi$ is a polynomial on $\F^{n+1}$ in real coordinates. This
$\psi$ must be invariant with respect to the rotations of $\F$, i.e.
$\psi(a\lambda) = \psi(a)$ for all $\lambda\in\F$, $\abs{\lambda}=1$.
It is not unique but becomes unique if it is required to be
homogeneous (which is always possible) of minimal degree. The latter
is said to be the {\em degree} of $\phi$. The number $\deg\phi$ is an
even integer since $\psi(-a)=\psi(a)$.

\begin{example}
  For every $t\in \N$ and every $y\in\fpn$ the function
  $\phi_{2t;y}(x) = (x,y)^t$, $x\in\fpn$, is a polynomial
  function of degree $2t$.
\end{example}

Given $d\in2\N$, we denote by $\polfd$ the space of all polynomial
functions of degrees $\leq d$. It has been proven in
\cite{lyushatpolyfunc} that the family $\{\phi_{d;y} :
y\in\fpn\}$ spans the whole space $\polfd$. We apply this
result to prove the following
\begin{prop}\label{prop:2.2}
  A finite nonempty set $X\subset\fpn$ is a tight $t$-design if and
  only if
  \begin{equation}
    \label{eq:2.2tight}
    \frac{1}{\abs{X}} \sum_{x\in X} \phi(x) = \int_{S^{2N-1}}
    \tilde\phi(a)\dm\sigma(a), \eqskip \phi\in\Pol_\F(2t),
  \end{equation}
  where $\tilde\phi$ is induced by the natural mapping
  $S^{2N-1}\rightarrow\fpn$ and $\sigma$ is the normalized Lebesgue measure.
\end{prop}

\begin{proof}
  The identity \eqref{eq:2.2tight} is equivalent to
  \begin{equation}
    \label{eq:2.3}
    \frac{1}{\abs{X}} \sum_{x\in X} F((x,y)) = \int_{S^{2N-1}} F(\abs{a^*b}^2)\dm\sigma(a),
  \end{equation}
  where $y=\hat b$, $b\in S^{2N-1}$, $F$ runs over the space $\Pi_t$
  of all univariate polynomials of degrees $\leq t$. By a known
  integration formula (see \cite{hoggardesignsproj}, Theorem 2.11) one
  can rewrite \eqref{eq:2.3} in the form
  \begin{equation}
    \label{eq:2.4}
    \frac{1}{\abs{X}} \sum_{x\in X} F((x,y)) = 
    \int_{-1}^1 F\left(\frac{1+\tau}{2}\right) \Omega_{\alpha,\beta}(\tau) \dm
    \tau, \eqskip \psi\in\Pi_t,
  \end{equation}
  where $\Omega_{\alpha,\beta}(\tau)$ is the normalized Jacobi weight,
  i.e.
  \begin{equation}
    \label{eq:2.4p}
    \Omega_{\alpha,\beta}(\tau) = c_{\alpha,\beta}(1-\tau)^\alpha
    (1+\tau)^\beta, \eqskip -1<\tau<1,
  \end{equation}
  with
  \begin{equation}
  \label{eq:3.10}
  c_{\alpha,\beta} = \left( \int_{-1}^1 (1-\tau)^\alpha
    (1+\tau)^\beta\dm\tau \right)^{-1} =
  \frac{\Gamma(\alpha+\beta+2)}{2^{\alpha+\beta+1}\Gamma(\alpha+1)\Gamma(\beta+1)}.
  \end{equation}
  In turn,
  \eqref{eq:2.4} is equivalent to its restriction to $F=P_i(\xi)$,
  $1\leq i\leq t$, since these polynomials constitute a basis in
  $\Pi_t$. It remains to note that
  \begin{equation*}
    \int_{-1}^1
    P_i\left(\frac{1+\tau}{2}\right)\Omega_{\alpha,\beta}(\tau) \dm \tau = 
    \int_{-1}^1 P_i^{(\alpha,\beta)}(\tau)\Omega_{\alpha,\beta}(\tau)\dm \tau = 0
  \end{equation*}
  by \eqref{eq:2.0pp} and the $\Omega_{\alpha,\beta}$-orthogonality of the
  system $\left\{ P_i^{(\alpha,\beta)} : i\geq 0 \right\}$.
\end{proof}

\begin{cor}\label{cor:2.3}
  Let $X\subset\fpn$ be a $t$-design. Then
  \begin{equation}
    \label{eq:2.5}
    \frac{1}{\abs{X}} \sum_{x\in X} P((u,x))Q((x,v)) = \int_{S^{2N-1}}
    P(\abs{a^*c}^2)Q(\abs{c^*b}^2)\dm\sigma(c)
  \end{equation}
  for $u=\hat a$, $v=\hat b$ and all univariate polynomials $P$, $Q$
  such that $\deg P +\deg Q\leq t$.
\end{cor}

\begin{proof}
  The mapping $x\mapsto P((u,x))Q((x,v))$, $x\in\fpn$, is a polynomial
  function of degree $\leq 2t$.
\end{proof}

\begin{cor}\label{cor:2.4}
  Let $X$, $P$, $Q$ be fixed under the conditions of Corollary
  \ref{cor:2.3}. Then the value
  \[
    \sum_{x\in X} P((u,x))Q((x,v))
  \]
  depends only on the inner product $(u,v)$ of $u,v\in\fpn$.
\end{cor}

\begin{proof}
  Let $(u_1,v_1)=(u,v)$, i.e. $\abs{a_1^*b_1}^2 = \abs{a^*b}^2$ where
  $\hat a_1=u_1$, $\hat b_1=v_1$. Without loss of generality one can
  assume that $a_1^*b_1 = a^*b$. Then there exists a $(n+1)\times(n+1)$
  matrix $T$ over $\F$ such that $T^*T=\id$ and $a_1=Ta$,
  $b_1=Tb$. This substitution in \eqref{eq:2.5} is equivalent to the
  change of variable $c\mapsto T^* c$. The latter does not affect the integral since
  the measure $\sigma$ is orthogonally invariant.
\end{proof}

\section{Bose-Mesner algebra}
\label{sec:3bosemesner}

Let $X$ be a finite nonempty subset of $\fpn$ and let
\[
  A'(X) = A(X)\cup\{1\} = \left\{(x,y) : x,y\in X\right\},
\]
so that $\abs{A'(X)} = s+1$. The $X\times X$ matrices of the form
\begin{equation}
  \label{eq:3.1}
  M_F = \left[F((x,y))\right]_{x,y\in X}
\end{equation}
where $F$ runs over all functions $A'(X)\rightarrow\C$, constitute a
complex linear space $\cD(X)$. Its natural basis consists of the
matrices
\begin{equation}
  \label{eq:3.1p}
  \Delta_\zeta = \left[ \delta_{\zeta,(x,y)}\right]_{x,y\in X},\eqskip
  \zeta\in A'(X),
\end{equation}
thus, $\dim\cD(X) = s+1$. The Lagrange interpolation formula allows us
to let $F$ in \eqref{eq:3.1} run over the polynomial space
$\Pi_s$, so that we have the isomorphism $F\mapsto M_F$ between
$\Pi_s$ and $\cD(X)$. In particular, if $F|A(X)=0$ and $F(1)=1$ then
$M_F=I$, the unit matrix.

According to Corollary \ref{cor:2.4} for $P,Q\in\Pi_s$, the matrix
product $M_P M_Q$ belongs to $\cD(X)$ if $\deg P + \deg Q\leq
t$. However, this condition is not fulfilled if $t=2s-1$ and $\deg P =
\deg Q = s$. Moreover, {\em Corollary \ref{cor:2.4} cannot be extended
  to this situation if $X$ is tight}. Indeed, suppose to the contrary
that
\[
  \sum_{x\in X} (u,x)^s(x,v)^s = \Phi((u,v)) \eqskip (u,v\in\fpn)
\]
with a function $\Phi:[0,1]\rightarrow \R_+$. Setting $v=u$ we obtain
\[
  \sum_{x\in X} (u,x)^{2s} = \Phi(1), \eqskip u\in\fpn.
\]
In other words,
\[
  \sum_{c\in\tilde X} \abs{a^* c}^{4s} = \Phi(1), \eqskip a\in
  S^{2N-1},
\]
where $\tilde X\subset S^{2N-1}$ is a complete system of
representatives of points $x\in X$, $\abs{\tilde X}=\abs{X}$. By integration
over $a$ we obtain
\[
  \Phi(1)= \left( \int_{S^{2N-1}} \abs{a^*c}^{4s}\dm\sigma(a)\right)
  \cdot \abs{X}
\]
since the integral does not depend on $c$. As a result,
\begin{equation}
  \label{eq:3.1pp}
  \frac{1}{\abs{X}} \sum_{x\in X} \phi_{4s;u}(x) = \int_{S^{2N-1}}
  \tilde\phi_{4s;u}(a)\dm\sigma(a),
\end{equation}
and by linearity, \eqref{eq:3.1pp} extends to the whole space
$\Pol_\F(4s)$. Thus, $X$ is a $2s$-design which is a contradiction
since $2s=t+1$.

Nevertheless,  under the constraint $u,v\in
X$, one can extend Corollary \ref{cor:2.4} to $t=2s-1$ and
$P,Q$ such that $\max(\deg P,\deg Q) = s$. This follows from the construction of a basis in $\cD(X)$ using
the Jacobi polynomials (cf. \cite{delsartegoethalsseidel}, Remark
7.6).

\begin{lemma}\label{lem:3.1}
  Let $X$ be a $t$-design in $\fpn$ and let
  $s=\left[\frac{t+1}{2}\right]$. Then $s+1$ matrices $M_i=M_{P_i}$,
  $0\leq i\leq s$, constitute a basis $\cM$ of $\cD(X)$ such that
  \begin{equation}
    \label{eq:3.2}
    M_i M_k =\abs{X} M_i\delta_{ik}\rho_{\mu(i,k)}
  \end{equation}
  where $\mu(i,k) = \min(i,k)$ and all $\rho_j>0$.
\end{lemma}

\begin{proof}
  The matrices $M_i$ are linearly independent because of the linear
  independence of the polynomials $P_i$. Since $\abs{\cM} = s+1$, this is a
  basis of $\cD(X)$. Now note that
  \[
    \int_{S^{2N-1}} P_i(\abs{a^*c}^2)P_k(\abs{c^*b}^2)\dm\sigma(c) =
    0, \eqskip i\neq k,
  \]
  by the addition formula for polynomial functions
  \cite{lyushatquaternions} (cf. \cite{hoggarpreprint},
  \cite{koornwinder}, \cite{muller}). The same formula with $i=k$ 
  yields
  \begin{equation}
    \label{eq:3.3}
    \int_{S^{2N-1}} P_i(\abs{a^*c}^2)P_i(\abs{c^*b}^2)\dm\sigma(c) =
    \chi_i P_i (\abs{a^* b}^2)
  \end{equation}
  where $\chi_i>0$. Assuming $\mu(i,k)\leq s-1$ (a fortiori,
  $i+k\leq 2s-1\leq t$) and using Corollary \ref{cor:2.3} we get
  \eqref{eq:3.2} with
  \begin{equation}
    \label{eq:3.4}
    \rho_j = \chi_j, \eqskip 0\leq j \leq s-1.
  \end{equation}
  In particular, $M_i M_s = M_s M_i =0$ for $0\leq i \leq s-1$. It
  remains to consider the case $i=k=s$.

  If $t$ is even the $t=2s$ and Corollary \ref{cor:2.3} is applicable
  to $i=k=s$, so $M^2_s = \abs{X}M_s\rho_s$ with 
  \begin{equation}
    \label{eq:3.5}
    \rho_s = \chi_s.
  \end{equation}
  
  Let $t$ be odd, so $t = 2s-1$. Then we decompose the unity matrix
  $I$ for the basis $\cM$,
  \begin{equation}
    \label{eq:3.5p}
    I = \sum_{i=0}^s \lambda_i M_i,
  \end{equation}
  and get $M_s = \lambda_s M_s^2$ multiplying \eqref{eq:3.5p} by
  $M_s$. This yields
  \begin{equation}
    \label{eq:3.5pp}
    \lambda_s = \frac{\tr M_s}{\tr M_s^2} =
    \frac{\abs{X}P_s(1)}{\sum_{x,y}P_s^2((x,y))} > 0,
  \end{equation}
  and then $M_s^2 =\abs{X}M_s\rho_s$ with
  \begin{equation}
    \label{eq:3.6}
    \rho_s = (\lambda_s\abs{X})^{-1}.
  \end{equation}
\end{proof}

\begin{remark}
  The formulas \eqref{eq:3.4} and \eqref{eq:3.5} are joined in
  \begin{equation}
    \label{eq:3.7}
    \rho_i = \chi_i, \eqskip 0\leq i \leq [t/2],
  \end{equation}
  while \eqref{eq:3.6} appears only for $t=2s-1$ in addition to
  \eqref{eq:3.7}.
\end{remark}

\begin{remark}
  The multiplication table \eqref{eq:3.2} shows that under conditions
  of Lemma \ref{lem:3.1} $\cD(X)$ is a commutative matrix algebra, the
  {\em Bose-Mesner algebra of $X$} \cite{bosemesner},
  \cite{delsartealgapproach}, \cite{delsartegoethalsseidel}.
\end{remark}

In what follows the conditions of Lemma \ref{lem:3.1} are
assumed to
be fulfilled. By setting
\begin{equation}
  \label{eq:3.8}
  L_i = \frac{M_i}{\rho_i\abs{X}} =
  \frac{1}{\rho_i\abs{X}}\left[P_i((x,y))\right]_{x,y\in X}
\end{equation}
the basis $\cM$ turns into $\cL = \left\{L_i\right\}_0^s$ consisting
of idempotents ($L_i^2 = L_i$) which are pairwise orthogonal ($L_iL_k
= 0$ for $i\neq k$). It is important to calculate their ranks.

We have
\[
  \rank L_i = \tr L_i = \rho_i^{-1}P_i(1),\eqskip 0\leq i\leq s,
\]
hence,
\[
  \rank L_i = \chi_i^{-1} P_i(1) = \frac{P_i^2(1)}{\int_{S^{2N-1}}
    P_i^2(\abs{a^*c}^2)\dm \sigma(c)}, \eqskip 0\leq i\leq [t/2]
\]
by \eqref{eq:3.7} and \eqref{eq:3.3} for $a=b$. Finally,
\begin{equation}
  \label{eq:3.9}
  \rank L_i = \frac{ \left(P_i^{(\alpha,\beta)}(1)\right)^2 }{
    \int_{-1}^1 \left(P_i^{(\alpha,\beta)}(\tau)\right)^2
    \Omega_{(\alpha,\beta)}(\tau) \dm\tau }, \eqskip 0\leq i\leq[t/2].
\end{equation}
In particular, $\rank L_0=1$. In addition to \eqref{eq:3.9} we have to find $\rank L_s$ in
the case $t=2s-1$. Formula \eqref{eq:3.6} is not effective to this end
since $\lambda_s$ is unknown. Indeed, in \eqref{eq:3.5pp} we cannot
proceed to the formally corresponding integral in the
denominator. Instead of this, we return to the decomposition of unity
and express $\rank L_s$ through $\rank L_i$, $0\leq i\leq
s-1$. We have
\[
  I = \sum_{i=0}^s L_i,
\]
whence,
\[
 \rank L_s = \tr L_s = \abs{X} - \sum_{i=0}^{s-1}\tr L_i = \abs{X} -
 \sum_{i=0}^{s-1} \rank L_i .
\]
By substitution from \eqref{eq:3.9} the last sum can be written as
$c_{\alpha,\beta}^{-1} K_{s-1}^{(\alpha,\beta)}(1,1)$, where
$K_{s-1}^{(\alpha,\beta)}(\cdot,\cdot)$ is the reproducing kernel of
the Jacobi polynomials with respect to the weight
$(1-\tau)^\alpha(1+\tau)^\beta$, see \cite{szegogabor}, Section
4.5. According to \eqref{eq:3.10} and formula (4.5.8) from
\cite{szegogabor} we obtain
\[
  \rank L_s = \abs{X} -
  \frac{\Gamma(s+\alpha+\beta+1)\Gamma(s+\alpha+1)\Gamma(\beta+1)}
  {\Gamma(\alpha+\beta+2)\Gamma(\alpha+2)\Gamma(s)}
\]
With our $\alpha,\beta$ defined by \eqref{eq:2.0palphabeta}
\begin{equation}
  \label{eq:3.11}
  \rank L_s = \abs{X} - \frac{\Gamma(N+s-1)\Gamma(N-m+s)\Gamma(m)}
  {\Gamma(N)\Gamma(N-m+1)\Gamma(m+s-1)\Gamma(s)} = \abs{X} -
  \frac{(N)_{s-1}(N-m+1)_{s-1}} {(m)_{s-1}(s-1)!}
\end{equation}

\begin{lemma}
  Let $X$ be a tight $t$-design in $\fpn$ with $t=2s-1$. Then
  \begin{equation}
    \label{eq:3.12}
    \rank L_s = \frac{(N)_{s-1}(N-m)_s}{(m)_s(s-1)!}
  \end{equation}
\end{lemma}
\begin{proof}
  In this case $e=s-1$, $\epsilon=1$, so \eqref{eq:2.1p}
  yields
  \begin{equation*}
    \abs{X} = R^1_{s-1}(1) = \frac{(N)_s(N-m+1)_{s-1}}{(m)_s (s-1)!}
  \end{equation*}
\end{proof}

The ranks of the other $L_i$ (including $L_s$ if $t=2s$) can be
explicitly calculated by \eqref{eq:3.9}, \eqref{eq:3.10} and \eqref{eq:2.0jacpoly} combined with (4.33) of
\cite{szegogabor}. This results in

\begin{lemma}
  Let $X$ be a $t$-design in $\fpn$ with
  $s=\left[\frac{t+1}{2}\right]$. Then
  \begin{equation}
    \label{eq:3.13}
    \rank L_i = \frac{ (N)_{i-1}(N-m)_i(N+2i-1)}{(m)_i i!}, \eqskip
    0\leq i\leq [t/2].
  \end{equation}
\end{lemma}

\begin{remark}
  Formula \eqref{eq:3.13} yields the true value $\rank L_0 =1$ by
  setting $(\gamma-1)(\gamma)_{-1} = 1$ for all $\gamma$.
\end{remark}

\begin{cor}
  The inequality
  \begin{equation}
    \label{eq:3.14}
    \rank L_i > \rank L_{i-1}, \eqskip 1\leq i \leq[t/2],
  \end{equation}
  holds, except for $X\subset\fp^1$. In the latter case
  \begin{equation}
    \label{eq:3.14p}
    \rank L_i =2, \eqskip 1\leq i \leq[t/2].
  \end{equation}
\end{cor}

Now note that our idempotents $L_i$ coincide with the matrices $E_i$
from \cite{hoggaroctonions} for $0\leq i\leq[t/2]$ but $L_s\neq E_s$
if $t=2s-1$, $X\not\subset\fp^1$. Indeed, according to (2.5)
from \cite{hoggaroctonions},
\begin{equation}
  \label{eq:3.15}
  E_i((x,y)) = \frac{1}{\abs{X}} \left[Q_i((x,y))\right]_{x,y\in X},
  \eqskip 0\leq i\leq s,
\end{equation}
where $Q_i(\xi)$ is proportional to $P_i(\xi)$ and
\begin{equation}
  \label{eq:3.16}
  Q_i(1) = \frac{(N)_{i-1} (N-m)_i (N+2i-1)} {(m)_i i!}, \eqskip i\geq
  0.
\end{equation}
Hence, $E_i$ are proportional to $L_i$ for all $i$, $0\leq i\leq
s$. Moreover, if $0\leq i\leq [t/2]$ then $\tr E_i = Q_i(1) = \tr L_i$
by \eqref{eq:3.16} and \eqref{eq:3.13}. Hence, $E_i = L_i$ for
$0\leq i \leq [t/2]$. However, if $t=2s-1$ (so $s=[t/2]+1$) and $X\not\subset\fp^1$
then $\tr E_s = Q_s(1) > \tr L_s$, see \eqref{eq:3.12}. In this case
$\tr E_s > \rank E_s$, so $E_s$  {\em is not an
  idempotent}. This is an obstacle to the full proof of
Theorem 2.6 \cite{hoggaroctonions} of the rationality of $A(X)$. To
overcome this difficulty, it suffices to change $E_s$ for $L_s$ (when
$t=2s-1$, $s\geq 2$) but then the ``critical inequality'' $\rank
L_s\neq \rank L_1$ is needed. However, the latter is not always
true. We clarify this intricate situation in the next section.

\section{The critical inequality and rationality theorem}
\label{sec:4critineq}

We prove the following
\begin{theorem}
  \label{thm:4.1}
  With $t=2s-1$, $s\geq 2$, the inequality
  \begin{equation}
    \label{eq:4.1}
    \rank L_s \neq \rank L_1
  \end{equation}
  holds for every tight $t$-design $X\subset\fpn$, except for a tight
  5-design in $\C\Proj^1$.
\end{theorem}
\begin{proof}
  From \eqref{eq:3.12} it follows that
  \[
    \rank L_s \geq \frac{(N)_1 (N-m)_2}{(m)_2\cdot 1!} =
    \frac{N(N-m)(N-m+1)}{m(m+1)}
  \]
  since the right side of \eqref{eq:3.12} increases with $s$. On the
  other hand, \eqref{eq:3.13} yields
  \begin{equation}
    \label{eq:4.2}
    \rank L_1 = \frac{ (N)_0(N-m)_1 (N+1)}{(m)_1\cdot 1!} = \frac{(N-m)(N+1)}{m}
  \end{equation}
  Hence,
  \[
    \rank L_s - \rank L_1 \geq \frac{(N-m)((N-m)^2 -
      (m^2+m+1))}{m(m+1)}
  \]
  Since $N-m=mn$ we obtain $\rank L_s > \rank L_1$ for $n^2 >
  1+m^{-1} + m^{-2}$, i.e. for $n\geq 3$ if $\F=\R$ and for $n\geq 2$
  if $\F=\C$ or $\Quat$. It remains to consider two cases. 
  \begin{enumerate}
    \item $\F = \R$, $n=2$. Then $m=1/2$, $N=3/2$, so $\rank L_1 = 5$
      by \eqref{eq:4.2}, while $\rank L_s = 2s$ by  \eqref{eq:3.12}. 
    \item $\F$ is arbitrary, $n=1$. Then $N=2m$, hence
  \[
    \rank L_s -\rank L_1 = \frac{(2m)_{s-1}} {(s-1)!} - (2m+1) =
    \begin{cases}
      -1 & \F = \R \\
      s-3& \F = \C \\
      \frac{1}{6}s(s+1)(s+2)-5 & \F=\Quat
    \end{cases}
  \]
   \end{enumerate}

    We see that $\rank L_s\neq\rank L_1$ with the only exception
  $\F=\C$, $s=3$, so $t=5$.
\end{proof}

A tight 5-design in $\C\Proj^1$ should contain 12 points since its
parameters are $m=1$, $N=2$, $s=3$, $e=2$, $\epsilon=1$. Accordingly,
\eqref{eq:2.1p} becomes
\[
  R_2^1(1) = \frac{(2)_3(2)_2}{(1)_3\cdot 2!} = 12.
\]
Such a design has been constructed in \cite{lyushatdedicata} as the
projective image of an orbit of the binary icosahedral group that is a
subgroup of $SU(2)$. Its representatives on the unit sphere
$S^3\subset\C^2$ are
\[
  a_1 = \left(\begin{matrix} 1 \\ 0\end{matrix}\right), a_2 =
  \left(\begin{matrix} 0 \\ 1\end{matrix}\right), a_k = \begin{cases}
    \mu \left(\begin{matrix} \lambda\eta^{k-3} \\ 1
      \end{matrix}\right) & 3\leq k \leq 7,\\ \vspace{.3cm}
    \mu\left(\begin{matrix} \eta^{k-3} \\ -\lambda\end{matrix}\right)
    & 8\leq k\leq 12, \end{cases}
\]
where
\[
  \eta = \exp\frac{2\pi i}{5},\eqskip \lambda=2\cos\frac{2\pi}{5} =
  \frac{\sqrt{5}-1}{2}, \eqskip \mu = \frac{1}{2\sin\frac{4\pi}{5}} =
  \sqrt{\frac{5+\sqrt{5}}{10}} .
\]
We omit an elementary calculation of the inner products $(x_j, x_k) =
\abs{a^*_j a_k}^2$, only noting that
\[
  \lambda^2 + \lambda -1 =0,\eqskip (2-\lambda)\mu^2 =1,\eqskip
  \abs{\eta^r -1}^2 = 0,2-\lambda, 3+\lambda,
\]
the latter for $r\equiv 0, \pm 1, \pm 2\mod 5$,
accordingly. A calculation yields,
\begin{equation}
  \label{eq:4.3}
  A(X) = \left\{0, \frac{5-\sqrt 5}{10}, \frac{5+\sqrt 5}{10}\right\},
\end{equation}
so $s=3$, $e=2$, $\epsilon=1$. The polynomial
\[
  \xi R^1_2(\xi) = 4\xi P_2^{(1,1)}(2\xi-1) = 6\xi(5\xi^2 - 5\xi +1 )
\]
annihilates $A(X)$. By Theorem C of Section 2 our $X$ is indeed a tight $5$-design.

The angle set \eqref{eq:4.3} is not rational. This occurs
because of the equality $\rank L_3 = \rank L_1$. In fact, 
$\rank L_0 = 1$, $\rank L_1 =3$, $\rank L_2 = 5$, $\rank L_3 =3$, by
\eqref{eq:3.12} and \eqref{eq:3.13}.

In the following corrected form of Theorem 2.6 \cite{hoggaroctonions}
the case $X\subset\R\Proj^1$ is also included for completeness. In this case, \eqref{eq:3.14p}
makes it possible for $A(X)$ to be not rational, though $\rank
L_s\neq\rank L_1$ under the conditions of Theorem \ref{thm:4.1}.

\begin{theorem}
  \label{thm:4.2}
  Let X be a tight $t$-design in $\fpn$. Then the angle set $A(X)$ is
  rational, except for two cases: 1) $X\subset\C\Proj^1$, $t=5$; 2)
  $X\subset\R\Proj^1$, $t\neq1,2,3,5$.
\end{theorem}
\begin{proof}
  For $X\not\subset\R\Proj^1$ the proof is the same as in
  \cite{hoggaroctonions} but with $L_s$ instead of $E_s$ when
  $t=2s-1$, $s\geq 2$, and using our Theorem \ref{thm:4.1} in this
  case.

  Now let $X\subset\R\Proj ^1$. Then it is the projective image of a
  regular $(2t+2)$-gon as easily follows from \cite{hong}. Therefore,
  \[
    A(X)\setminus\{0\} = \left\{ \cos^2\frac{k\pi}{t+1}\right\}_1^e =
    \left\{\frac{1}{2}
      \left(1+\cos\frac{2k\pi}{t+1}\right)\right\}_1^e
  \]
  where $e=[t/2]$. Since $\cos m\theta$ is a polynomial of $\cos
  \theta$ with integer coefficients, the set $A(X)$ is rational if and
  only if the number $\rho =\cos\frac{2\pi}{t+1}$ is
  rational. Obviously, the latter is true if $t=1,2,3,5$. Conversely,
  let $\rho\in\Q$. Then the complex number $w=\exp(2\pi i/(t+1))$
  satisfies the equation $w^2 -2\rho w + 1=0$. On the other hand, this
  is a primitive root of 1 of degree $t+1$. It is known that
  irreducible (over $\Q$) equation for $w$ is of degree $\varphi(t+1)$
  where $\varphi$ is the Euler function. Hence, $\varphi(t+1)\leq2$. If
  $\varphi(t+1)=1$ then $t=1$. If $\varphi(t+1)=2$ then $t\in \{2,3,5\}$ as
  easily follows from the classical formula
  \[
     \varphi\left(\prod_{i=1}^r q_i^{\nu_i}\right) = \prod_{i=1}^r
     q_i^{\nu_i -1}(q_i -1)
  \]
  where $q_1,\ldots,q_r$ are prime divisors of $t+1$.
\end{proof}

\nocite{*}
\bibliographystyle{plain}
\bibliography{tightdesigns}

Department of Mathematics

Technion, Haifa 32000, Israel

email: lyubich@tx.technion.ac.il

\end{document}